\RequirePackage[l2tabu, orthodox]{nag}
\documentclass[10pt, draft]{article}
\usepackage{dsfont}
\usepackage{amsmath,amssymb,amsbsy,amsfonts,amsthm,latexsym,amsopn,amstext,
                            amsxtra,euscript,amscd}
\usepackage{booktabs}
\usepackage{enumitem}
\usepackage{tikz}
\usepackage[norefs,nocites]{refcheck}

\newtheorem{thm}{Theorem}

\def\\{\cr}
\def\({\left(}
\def\){\right)}
\def\[{\left[}
\def\]{\right]}
\def\<{\langle}
\def\>{\rangle}
\def\fl#1{\left\lfloor#1\right\rfloor}

\def\notdivides{\mathrel{\kern-3pt\not\!\kern3.5pt\bigm|}}

\begin{document}


\title{\Large \textbf{On Some Floor Function Sets}}
\author{
\scshape {RANDELL HEYMAN}  \\
School of Mathematics and Statistics,\\ University of New South Wales,\\
Sydney, Australia\\
\texttt {randell@unsw.edu.au}
\and
\scshape{MD RAHIL MIRAJ}\\
{Theoretical Statistics and Mathematics Unit},\\
 {Indian Statistical Institute}, \\
{Bengaluru, Karnataka, India}\\
\texttt{ rahilmiraj@gmail.com}
}
\date{\today}

\maketitle



\begin{abstract}

Let $X$ be a positive integer and $t$ a real number great than 1. The family of sets
$\left\{\fl{\frac{X}{n^t}} ~:~ 1\leq n\leq X\right\}$ have an interesting prime distribution property. We give an exact formula for the cardinality of these sets. We provide an estimate for the cardinality of the set $\left\{\fl{\frac{X}{p}} ~:~ p~ \text{prime},~ p\leq X\right\}$. For positive real $X$, we derive asymptotic formulas for the cardinality of the set $\big\{\fl{f(n)} ~:~ 1\leq n\leq X\big\}$ for various sets of functions.

\end{abstract}
\newline
AMS Classification (2020): 11A41, 11A67, 11B05

\section{Introduction}

Sets and sequences utilising the floor function, and in particular the primality of their elements, have been studied for many decades. The two most commonly researched are probably Beatty sequences (see, for example, ~\cite{ABS,BaBa,BaLi,GuNe,Harm}) and Piatetski-Shapiro sequences (see, for example,~\cite{Akb,BBBSW,BBGY,BGS,LSZ,Morg,RivW}).
Recent research interest has also focused on `hyperbolic sets' utilising the floor function.

For positive integer $X$ let
$$S(X):= \left\{\left\lfloor \frac{X}{n} \right\rfloor ~:~ 1\leq n\leq X\right\}.$$
In \cite{Hey} it was shown that
$$|S(X)|=2\sqrt{X}+O(1).$$
An exact formula was also given for integer values of $X$.
Specifically,
$$|S(X)|= \fl{\sqrt{4X+1}}-1.$$
Using exponential sum techniques, it was shown (in \cite{Hey2}) that the number of primes in the set $S(X)$ could be estimated. That is,
$$|\{p \in S(X): p \text{ is prime}\}|=\frac{4\sqrt{X}}{\log X}+O\(\frac{\sqrt{X}}{(\log X)^2}\),$$
which was improved by Ma And Wu \cite{Ma3} as follows:
\begin{align*}
|\{p \in S(X): p \text{ is prime}\}|
 &=\int_2^{\sqrt{X}}\frac{dt}{\log t}+\int_2^{\sqrt{X}}\frac{dt}{\log (X/t)}\\&+O\(\sqrt{X}\exp\(-c(\log X)^{3/5}(\log \log X)^{-1/5}\)\),
\end{align*}
where $c>0$ is a positive integer.
In that paper Ma And Wu observed that $S(X)$ is probably the first example of such a sparse set that satisfies the prime number theorem, in the sense that
$$|\{p \in S(X): p \text{ is prime}\}|\sim\frac{|S(X)|}{\log |S(X)| }.$$
More recently, in \cite{Bor}, it was shown that there exists even sparser sets that satisfy the prime number theorem. For real $t>1$ the family of sets
$$S_t(X):= \left\{\left\lfloor \frac{X}{n^t} \right\rfloor ~:~ 1\leq n\leq X\right\},$$
also satisfy the prime number theorem.

In this paper we investigate some variations that naturally arise from the sets $S(X)$ and $S_t(X)$.

We start with the cardinality of $S_t(X)$. An estimate of the cardinality of $S(X)$ is given in \cite{Bor} as follows.
$$|S_t(X)|=X^{\frac1{t+1}}\(t^{\frac{t}{t+1}}+t^{\frac1{t+1}}\)+O(1).$$
Our first theorem gives an exact formula.
\begin{thm}
\label{thm:exact}
For given integer $X\geq 1$ and $t>1$, let $a=(tX)^{1/(t+1)}$. If $X=1$ we have $|S_t(X)|=1$. For $X\geq 2$ we have
$$|S_t(X)|=
\begin{cases}
 a+\left\lfloor{\frac{X}{a^t}}\right\rfloor   & \text{if $a$ is integer},\\
 \lfloor{a}\rfloor+\left\lfloor{\frac{X}{\lfloor{a+1}\rfloor^t}}\right\rfloor+\varepsilon(X,t)
 & \text{if $a$ is not an integer},\\
\end{cases}$$
where
$$\varepsilon(X,t)=
\begin{cases}
0  & \text{if $X\ge \left\lfloor{\frac{X}{\lfloor{a}\rfloor^t}}\right\rfloor \lfloor{a+1}\rfloor^t$},\\
1  & \text{otherwise.}
\end{cases}$$
\end{thm}

Our next result restricts the $n$ in $S(X)$ to primes. This theorem is implied by a careful reading of \cite{Erd}, but we include a full proof in section
Theorem \ref{thm:prime}.
\begin{thm}
\label{thm:prime}
Let
$$S_p(X)=\left\{\fl{\frac{X}{p}}: p~\textrm{prime}, ~p \le X \right\}.$$
There exists positive reals $c_1$ and $c_2$ such that
$$c_1 \(\frac{X}{\log X}\)^{1/2}\le |S_p(X)|\le c_2 \(\frac{X}{\log X}\)^{1/2}.$$
\end{thm}


Our final result generalises from the `hyperbolic' functions imbedded in $S(X)$ and $S_t(X)$ to other functions.

\begin{thm}
\label{thm:general}
For an arbitrary positive real $X$ let
$$S_f (X)=\big\{\fl{f(n)}: 1\le n \le X\big\},$$
where $f(y)=f(X,y)$ is a non-negative twice differentiable function and there exists a unique value $a=a(X) \in[1,X]$ such that $\frac{d f}{d y}\mid_{y=a}=\pm1$. Then

i. if, for all $y \in [1,X]$, we have $f'(y)< 0$ and $f''(y)>0$ then
$$|S_f(X)|=f(\fl{a+1})-f(\fl{X})+a+O(1),$$

ii. if, for all $y \in [1,X]$, we have $f'(y)< 0$ and $f''(y)<0$ then
$$|S_f (X)|=f(1)-f(\fl{a})+X-a+O(1),$$

iii. if, for all $y \in [1,X]$, we have $f'(y)> 0$ and $f''(y)>0$ then
$$|S_f (X)|=f(\fl{a})-f(1)+X-a+O(1),$$

iv. if, for all $y \in [1,X]$, we have $f'(y)> 0$ and $f''(y)<0$ then
$$|S_f (X)|=f(\fl{X})-f(\fl{a+1})+a+O(1).$$
\end{thm}






As an example of the use of this theorem consider the `circle' set
$$\left\{\fl{\sqrt{X^2-n^2}}: 1\le n \le X\right\}.$$

Since for all $y \in [1,X]$ we have $f'(y)<0$ and $f''(y)<0$, we calculate the unique value $a=X/\sqrt{2}$. Then, by Theorem \ref{thm:general}, we obtain
\begin{align*}
\left|\left\{\fl{\sqrt{X^2-n^2}}: 1\le n \le X\right\}\right |
&=f(1)-f(\fl{a})+X-a+O(1)\\
&=(2-\sqrt{2})X+O(1).
\end{align*}

\section{Proof of Theorem \ref{thm:exact}}

It is easy to see that for $X=1$ we have $S_t(X)=\{1\}$ for any $t>1$.

We base our proof on the graph of $f(n)=X/n^t$, noting that values $\fl{X/n^t}$ can be inferred from the graph.
Let $a$ be the (unique) value of $n$ for which $f'(n)=-1$. Simple calculus calculations show that
$a=(tX)^{1/(t+1)}$.

Consider the case when $a$ is an integer. For $1 \le n < a$, utilising the mean value theorem, we see that
$$\frac{X}{(n+1)^t}-\frac{X}{n^t}<-1$$ and thus
$$\fl{\frac{X}{(n+1)^t}}-\fl{\frac{X}{n^t}}<\frac{X}{(n+1)^t}-\frac{X}{n^t}+1<0.$$
It follows that $\fl{\frac{X}{n^t}}$ and $\fl{\frac{X}{(n+1)^t}}$ are distinct. Therefore
$$\fl{\frac{X}{1^t}},\fl{\frac{X}{2^t}},\ldots,\fl{\frac{X}{(a-1)^t}} \in S_t(X).$$
So, for $1 \le n < a$, the contribution to the cardinality of $S_t(X)$ is $a-1$.

For $a \le n \le X$, the maximum element of $S_t(X)$ is $\fl{\frac{X}{a^t}}$ and the minimum element is $0$. So in this case we obtain
$$\fl{\frac{X}{a^t}}, \fl{\frac{X}{a^t}}-1,\ldots,1,0 \in S_t(X).$$
So, for $a \le n \le X$, the contribution to the cardinality of $S_t(X)$ is $\fl{\frac{X}{a^t}}+1.$

Again using the mean value theorem, $\fl{\frac{X}{a^t}}\ne \fl{\frac{X}{(a-1)^t}}$. So there is no over counting. More precisely,
$$\left\{\fl{\frac{X}{1^t}},\fl{\frac{X}{2^t}},\ldots,\fl{\frac{X}{(a-1)^t}}\right\} \bigcap \left \{\fl{\frac{X}{a^t}}, \fl{\frac{X}{a^t}}-1,\ldots,1,0\right \} = \phi$$
We conclude that if $a$ is an integer then
$$|S_t(X)|=a+\fl{\frac{X}{a^t}}.$$

Next, we examine the other case, $a$ is not an integer. Using a similar process we see that for $1 \le n \le \fl{a}$ we have $f'(n) <  -1$. So
$$\fl{\frac{X}{1^t}},\fl{\frac{X}{2^t}}, \ldots , \fl{\frac{X}{\fl{a}^t}} \in S_t(X).$$
Therefore, the values of $n$ for which $1\le n \le \fl{a}$ contribute to the cardinality of $S_t(X)$ is $\fl{a}$.

For $\fl{a} < n \le X$ we have
$$\fl{\frac{X}{\fl{a+1}^t}},\fl{\frac{X}{\fl{a+1}^t}}-1, \ldots , 1,0\in S_t(X).$$
So the values of $n$ for which $\fl{a} < n \le X$ contribute to the cardinality of $S_t(X)$ is $\fl{\frac{X}{\fl{a+1}^t}}$+1.
But it maybe that
$$\fl{\frac{X}{\fl{a}^t}}=\fl{\frac{X}{\fl{a+1}^t}}.$$
So we correct for this possible overlap by  1 if
$$\fl{\frac{X}{\fl{a}^t}}=\fl{\frac{X}{\fl{a+1}^t}},$$
and making no change otherwise. This will require two divisions, which takes $O(n^2)$ time each.
We will show that 1 should be subtracted from the cardinality if, and only if, $X\ge \left\lfloor{\frac{X}{\lfloor{a}\rfloor^t}}\right\rfloor \lfloor{a+1}\rfloor^t$. Having to now perform one division and one multiplication, the checking for overlaps will be slightly faster than before, as multiplication takes $O(n^{1.58})$ time following Karatsuba's algorithm \cite{Kars}.

If 1 should be subtracted from the cardinality then an overlap exists and then, for some non-negative integer $k$, we have
\begin{align}
\label{eq:k}
\fl{\frac{X}{\fl{a}^t}}&=\fl{\frac{X}{\fl{a+1}^t}}=k.
\end{align}

Then
$$\fl{\frac{X}{\fl{a}^t}}=k=\fl{\frac{X}{\fl{a+1}^t}}\leq \frac{X}{\fl{a+1}^t},$$
and so
$$ \fl{\frac{X}{\fl{a}^t}}\leq \frac{X}{\fl{a+1}^t},$$
which implies that we should deduct 1 if
$$X\ge \left\lfloor{\frac{X}{\lfloor{a}\rfloor^t}}\right\rfloor \lfloor{a+1}\rfloor^t.$$
Conversely,
$$X\ge \left\lfloor{\frac{X}{\lfloor{a}\rfloor^t}}\right\rfloor \lfloor{a+1}\rfloor^t$$
 implies
$$ \fl{\frac{X}{\fl{a}^t}}\leq \frac{X}{\fl{a+1}^t},$$
which means that
$$\fl{\frac{X}{\fl{a}^t}}\leq \frac{X}{\fl{a+1}^t}<\frac{X}{\fl{a}^t}<\fl{\frac{X}{\fl{a}^t}}+1.$$
From \eqref{eq:k} we have $k=\fl{\frac{X}{\fl{a}^t}}$, and therefore
$$k \le \frac{X}{\fl{a+1}^t}<\fl{\frac{X}{\fl{a}^t}}+1=k+1. $$
This means that \textit{}
$$k \le  \frac{X}{\fl{a+1}^t}< k+1,$$
which implies that
$$k=\fl{\frac{X}{\fl{a+1}^t}},$$
and so an overlap exists, and 1 should be subtracted from the cardinality.


We conclude that if $a$ is not an integer then

$$|S_t(X)|=\lfloor{a}\rfloor+\left\lfloor{\frac{X}{\lfloor{a+1}\rfloor^t}}\right\rfloor+\varepsilon(X,t),$$
where
$$\varepsilon(X,t)=
\begin{cases}
0  & \text{if $X\ge \left\lfloor{\frac{X}{\lfloor{a}\rfloor^t}}\right\rfloor \lfloor{a+1}\rfloor^t$},\\
1  & \text{otherwise,}
\end{cases}$$
concluding the proof.

\section{Proof of Theorem \ref{thm:prime}}
We denote by $v_p(M)$ the exponent of prime $p$ in the prime factorisation of a positive integer $M$.
For positive integer $M$ let
    $$\alpha(M)=\{\alpha_1, \ldots, \alpha_k:M=p_1 ^{\alpha_1}\cdots p_k ^{\alpha_k}\},$$
    where
    $M=p_1 ^{\alpha_1}\cdots p_k ^{\alpha_k}$ with  $p_1<p_2< \cdots < p_k$ is the prime factorisation of $M$  throughout.
We will show that
$$|S_p(X)|=|\alpha(X!)|.$$
Observe that this means that the cardinality of $S_p(X)$ is equal to the number of distinct exponents in the prime factorisation of $X!$.

Firstly, we show that the equality is true for the smaller primes. That is
\begin{align}
    \label{eq:small primes}
&\Bigg|\left\{\fl{\frac{X}{p}}: p \textrm{ prime }, p \le \sqrt{X} \right\}\Bigg|\notag\\
&=\Big|\left\{\alpha_1, \ldots, \alpha_j:X!=p_1 ^{\alpha_1}\cdots p_j^{\alpha_j} p_{j+1}^{\alpha_{j+1}} \cdots p_k ^{\alpha_k}, p_j \le \sqrt{X}< p_{j+1}\right\}\Big|.
\end{align}
For any two primes $2 \le p<q \le \sqrt{X}$ we have, using the methodology of Theorem \ref{thm:exact}, that
$$\fl{\frac{X}{p}}>\fl{\frac{X}{q}}.$$

So,
$$\Bigg|\left\{\fl{\frac{X}{p}}: p \textrm{ prime }, p \le \sqrt{X} \right\}\Bigg|=\pi(\sqrt{X}).$$
It is also clear that for $t \ge 2$ we have
$$\fl{\frac{X}{p^t}}\ge\fl{\frac{X}{q^t}}.$$
Next,
\begin{align*}
v_p(X!)=\sum_{t=1}^{\infty}\fl{\frac{X}{p^t}}>\sum_{t=1}^{\infty}\fl{\frac{X}{q^t}}=v_q(X!).
\end{align*}
This means that every prime not exceeding $\sqrt{X}$ adds exactly 1 to the cardinality of
$$\Big|\left\{\alpha_1, \ldots, \alpha_j:X!=p_1 ^{\alpha_1}\cdots p_j^{\alpha_j} p_{j+1}^{\alpha_{j+1}} \cdots p_k ^{\alpha_k}, p_j \le \sqrt{X}< p_{j+1}\right\}\Big|.$$
Therefore,
$$\Big|\left\{\alpha_1, \ldots, \alpha_j:X!=p_1 ^{\alpha_1}\cdots p_j^{\alpha_j} p_{j+1}^{\alpha_{j+1}} \cdots p_k ^{\alpha_k}, p_j \le \sqrt{X}< p_{j+1}\right\}\Big|=\pi(\sqrt{X}),$$
which proves equation \eqref{eq:small primes}.

Secondly, we show that $|S_p(X)|=|\alpha(X!)|$ is true with respect to the larger primes. Specifically,
\begin{align}
    \label{eq:large primes}
&\Bigg|\left\{\fl{\frac{X}{p}}: p \textrm{ prime }, \sqrt{X} <p \le X \right\}\Bigg|\notag\\
&=\Big|\left\{\alpha_{j+1}, \ldots, \alpha_k:X!=p_1 ^{\alpha_1}\cdots p_j^{\alpha_j} p_{j+1}^{\alpha_{j+1}} \cdots p_k ^{\alpha_k}, p_j \le \sqrt{X}< p_{j+1}\right\}\Big|.
\end{align}
For any prime $p > \sqrt{X}$ and $t>1$ we have $\frac{X}{p^t}<1$. So
$$v_p(X!)=\sum_{t=1}^{\infty}\fl{\frac{X}{p^t}}=\fl{\frac{X}{p}}.$$
Thus
\begin{align*}
    &\left\{\fl{\frac{X}{p}}: p \textrm{ prime }, \sqrt{X} <p \le X \right\}\notag\\
&=\left\{\alpha_{j+1}, \ldots, \alpha_k:X!=p_1 ^{\alpha_1}\cdots p_j^{\alpha_j} p_{j+1}^{\alpha_{j+1}} \cdots p_k ^{\alpha_k}, p_j \le \sqrt{X}< p_{j+1}\right\},
\end{align*}
which means that equality \eqref{eq:large primes} is true.
Having shown that $|S_p(X)|=|\alpha(X!)|$ we note (see \cite{Erd}) that the exists positive reals $c_1,c_2$ such that
$$c_1 \(\frac{X}{\log X}\)^{1/2}\le |\alpha(X!)| \le c_2\(\frac{X}{\log X}\)^{1/2},$$
which concludes the proof.

\section{Proof of Theorem \ref{thm:general}}
Recall that $a$ is the unique value for which $\frac{df}{dn}$ equals -1 or 1. We use the mean value theorem freely throughout this section.

For the case where $f'(y)<0$ and $f''(y)>0$, we have $f'(a)=-1$. First consider when $1<n \le \fl{a}$. For some $c \in (n-1,n)$ we have
$$\frac{f(n)-f(n-1)}{n-(n-1)}=f'(c)< f'(a)=-1.$$
Thus, since $x-1<\fl{x}\le x$, $f(n-1)-f(n)>1$ implies
$$\fl{f(n-1)}-\fl{f(n)}>f(n-1)-1-f(n)>0.$$
Therefore,
$$\fl{f(1)},\fl{f(2)},\ldots,\fl{f(\fl{a})}\in S_f(X).$$
This contributes $\fl{a}$ to the cardinality of $S_f(X)$.

Next, consider when $\fl{a+1}\le n< X$. For some $d \in (n,n+1)$ we have
$$\frac{f(n+1)-f(n)}{(n+1)-n}=f'(d)>f'(a)=-1.$$
Thus $f(n)-f(n+1)<1$, implying that
$$\fl{f(n)}-\fl{f(n+1)}<f(n)-f(n+1)+1<2.$$
Hence we obtain $\fl{f(n)}-\fl{f(n+1)}\in\{0,1\}$.
Therefore, every integer between $f(\fl{X})$ and $f(\fl{a+1})$ inclusive is an element of $S_f(X)$. This contributes $\fl{f(\fl{X})}-\fl{f(\fl{a+1})}$ to the cardinality of $S_f(X)$.

Finally, there may be an overlap of at most one element between the sets $\big\{\fl{f(n)}: 1\le n \le \fl{a}\big\}$ and $\big\{\fl{f(n)}: \fl{a+1}\le n \le X\big\}$, if $\fl{f(\fl{a})}$ = $\fl{f(\fl{a+1})}$.
Combining all the above proves
$$|S_f(X)|=f(\fl{a+1})-f(\fl{X})+a+O(1)$$
for this case.

Now, if $f'(y)<0$ and $f''(y)<0$, we have $f'(a)=-1$. We consider when $\fl{a+1}\le n< X$. For some $c \in (n,n+1)$ we have
$$\frac{f(n+1)-f(n)}{(n+1)-n}=f'(c)<f'(a)=-1.$$
Thus, $f(n)-f(n+1)>1$, implying that
$$\fl{f(n)}-\fl{f(n+1)}>f(n)-1-f(n+1)>0.$$
Therefore,
$$\fl{f(\fl{a+1})},\fl{f(\fl{a+2})},\ldots,\fl{f(\fl{X})}\in S_f(X).$$
This contributes $\fl{X-a}$ to the cardinality of $S_f(X)$.

Next, consider when $1<n \le \fl{a}$. For some $d \in (n-1,n)$ we have
$$\frac{f(n)-f(n-1)}{n-(n-1)}=f'(d)> f'(a)=-1.$$
Thus, $f(n-1)-f(n)<1$ implies
$$\fl{f(n-1)}-\fl{f(n)}<f(n-1)-f(n)+1<2.$$
Hence we obtain $\fl{f(n)}-\fl{f(n+1)}\in\{0,1\}$. Therefore, every integer between $f(1)$ and $f(\fl{a})$ inclusive is an element of $S_f(X)$. This contributes $\fl{f(1)}-\fl{f(\fl{a})}$ to the cardinality of $S_f(X)$.

Finally, there may be an overlap of at most one element between the sets $\big\{\fl{f(n)}: 1\le n \le \fl{a}\big\}$ and $\big\{\fl{f(n)}: \fl{a+1}\le n \le X\big\}$, if $\fl{f(\fl{a})}$ = $\fl{f(\fl{a+1})}$.
Combining all the above proves
$$|S_f(X)|=f(1)-f(\fl{a})+X-a+O(1)$$
for this case.

Next, we will consider the case when $f'(y)>0$ and $f''(y)>0$. We have $f'(a)=1$. We consider when $\fl{a+1}\le n< X$. For some $c \in (n,n+1)$ we have
$$\frac{f(n+1)-f(n)}{(n+1)-n}=f'(c)>f'(a)=1.$$
Thus, $f(n+1)-f(n)>1$, implying that
$$\fl{f(n+1)}-\fl{f(n)}>f(n+1)-1-f(n)>0.$$
Therefore,
$$\fl{f(\fl{a+1})},\fl{f(\fl{a+2})},\ldots,\fl{f(\fl{X})}\in S_f(X).$$
This contributes $\fl{X-a}$ to the cardinality of $S_f(X)$.

Next, consider when $1<n \le \fl{a}$. For some $d \in (n-1,n)$ we have
$$\frac{f(n)-f(n-1)}{n-(n-1)}=f'(d)< f'(a)=1.$$
Thus, $f(n)-f(n-1)<1$ implies
$$\fl{f(n)}-\fl{f(n-1)}<f(n)-f(n-1)+1<2.$$
Hence we obtain $\fl{f(n)}-\fl{f(n-1)}\in\{0,1\}$. Therefore, every integer between $f(1)$ and $f(\fl{a})$ inclusive is an element of $S_f(X)$. This contributes $\fl{f(1)}-\fl{f(\fl{a})}$ to the cardinality of $S_f(X)$.

Finally, there may be an overlap of at most one element between the sets $\big\{\fl{f(n)}: 1\le n \le \fl{a}\big\}$ and $\big\{\fl{f(n)}: \fl{a+1}\le n \le X\big\}$, if $\fl{f(\fl{a})}$ = $\fl{f(\fl{a+1})}$.
Combining all the above proves
$$|S_f(X)|=f(\fl{a})-f(1)+X-a+O(1)$$
for this case.

Our last case is when $f'(y)>0$ and $f''(y)<0$. We have $f'(a)=1$. We first consider when $1<n \le \fl{a}$. For some $c \in (n-1,n)$ we have
$$\frac{f(n)-f(n-1)}{n-(n-1)}=f'(c)> f'(a)=1.$$
Therefore, $f(n)-f(n-1)>1$ implies
$$\fl{f(n)}-\fl{f(n-1)}>f(n-1)-1-f(n)>0.$$
Hence,
$$\fl{f(1)},\fl{f(2)},\ldots,\fl{f(\fl{a})}\in S_f(X).$$
This contributes $\fl{a}$ to the cardinality of $S_f(X)$.

Next, consider when $\fl{a+1}\le n< X$. For some $d \in (n,n+1)$ we have
$$\frac{f(n+1)-f(n)}{(n+1)-n}=f'(d)<f'(a)=1.$$
Thus, $f(n+1)-f(n)<1$, implying that
$$\fl{f(n+1)}-\fl{f(n)}<f(n+1)-f(n)+1<2.$$
Hence we obtain $\fl{f(n+1)}-\fl{f(n)}\in\{0,1\}$. Therefore, every integer between $f(\fl{a+1})$ and $f(\fl{X})$ inclusive is an element of $S_f(X)$. This contributes $\fl{f(\fl{a+1})}-\fl{f(\fl{X})}$ to the cardinality of $S_f(X)$.

Finally, there may be an overlap of at most one element between the sets $\big\{\fl{f(n)}: 1\le n \le \fl{a}\big\}$ and $\big\{\fl{f(n)}: \fl{a+1}\le n \le X\big\}$, if $\fl{f(\fl{a})}$ = $\fl{f(\fl{a+1})}$.
Combining all the above proves
$$|S_f(X)|=f(\fl{X})-f(\fl{a+1})+a+O(1)$$
for this case.

\section{Acknowledgment}
We thank an anonymous referee for useful comments. Additionally, he or she supplied an exact formula for our final theorem.
\makeatletter
\renewcommand{\@biblabel}[1]{[#1]\hfill}
\makeatother

\end{document}